# ON MIXED TYPE DUALITY FOR NONDIFFERENTIABLE MULTIOBJECTIVE VARIATIONAL PROBLEMS

**I.Husain**

Department of Mathematics, Jaypee University of Engineering and Technology, Guna, MP, India

Email:ihusain11@yahoo.com.

**Rumana, G.Mattoo**

Department of Statistics, University of Kashmir, Srinagar, Kashmir, India.

## ABSTRACT

A mixed type dual to a nondifferentiable variational problem involving higher order derivative is formulated and duality results are proved under generalized invexity conditions. Special cases are generated from our results.

## 1. INTRODUCTION

In [9], Husain et al considered the following nondifferentiable multiobjective variational problem containing square root terms:

**(VP)** Minimize
$$\left( \int_I \left( f^1(t,x,\dot{x},\ddot{x}) + \left( x(t)^T B^1(t) x(t) \right)^{\frac{1}{2}} \right) dt, \ldots, \int_I \left( f^p(t,x,\dot{x},\ddot{x}) + \left( x(t)^T B^p(t) x(t) \right)^{\frac{1}{2}} \right) dt \right)$$

Subject to

$$x(a) = 0 = x(b)$$

$$\dot{x}(a) = 0 = \dot{x}(b)$$

$$g(t,x,\dot{x},\dot{x}) \leqq 0 \ , \ t \in I,$$

where

1) $f^i : I \times R^n \times R^n \times R^n \to R, (i = 1,2,\ldots,p)$ , $g : I \times R^n \times R^n \times R^n \to R^m$ , $j = 1,\ldots,m$ are assumed to be continuously differentiable functions, and for each $t \in I$, $i \in P = \{1,2,\ldots,p\}$ , $B^i(t)$ is an $n \times n$ positive semidefinite (symmetric) matrix, with $B^i(\cdot)$ continuous on $I$ .

The following Wolfe type dual to the problem (VP) is presented in [9]:

**(MWD)** Maximize$\Bigg( \int_I \Big( f^1(t,u,\dot{u},\ddot{u}) + \big(u(t)^T B^1(t) z^1(t)\big) + y(t)^T g(t,u,\dot{u},\ddot{u})\Big) dt$

$$,\ldots, \int_I \Big( f^p(t,u,\dot{u},\ddot{u}) + \big(u(t)^T B^p(t) z^p(t)\big) + y(t)^T g(t,u,\dot{u},\ddot{u})\Big) dt \Bigg)$$

Subject to

$$u(a) = 0 = u(b)$$

$$\dot{u}(a) = 0 = \dot{u}(b)$$

$$\sum_{i=1}^{p} \lambda^i \Big( f_u^i(t,u,\dot{u},\ddot{u}) dt + B^i(t) z^i(t) + y(t)^T g_u(t,u,\dot{u},\ddot{u})\Big)$$

$$- D\Big(\lambda^T f_{\dot{u}}(t,u,\dot{u},\ddot{u}) + y(t)^T g_{\dot{u}}(t,u,\dot{u},\ddot{u})\Big)$$

$$+ D^2\Big(\lambda^T f_{\ddot{u}}(t,u,\dot{u},\ddot{u}) + y(t)^T g_{\ddot{u}}(t,u,\dot{u},\ddot{u})\Big) = 0 \ , \ t \in I,$$

$$\bar{z}^i(t)^T B^i(t) \bar{z}^i(t) \leqq 1 \ , \ t \in I \ , \ i \in P,$$

$$y(t) \geqq 0 \ , \ t \in I \ ,$$

$$\lambda > 0 \ , \ \lambda^T e = 1.$$

The problem (MWD) is a dual to (VP) assuming that

$$\sum_{i=1}^{p} \lambda^i \int_I \Big( f^i(t,.,.,.) + (\cdot)^T B^i(t) z^i(t) + y(t)^T g(t,.,.,.)\Big) dt$$

is pseudoinvex with respect to $\eta$ . The authors in [9] further weakened the invexity required in Wolfe type duality by constructing the following Mond- Weir type vector dual.

The following is the Mond-Weir vector type dual to (VP):

**(M-WVD)** Maximize$\left( \int_I \left( f^1(t,u,\dot{u},\ddot{u}) + \left(u(t)^T B^1(t) z^1(t)\right)\right) dt \right.$

$$\left. ,\ldots, \int_I \left( f^p(t,u,\dot{u},\ddot{u}) + \left(u(t)^T B^p(t) z^p(t)\right)\right) dt \right)$$

Subject to

$$u(a) = 0 = u(b)$$

$$\dot{u}(a) = 0 = \dot{u}(b)$$

$$\sum_{i=1}^{p} \lambda^i \left( f_u^i(t,u,\dot{u},\ddot{u}) dt + B^i(t) z^i(t) + y(t)^T g_u(t,u,\dot{u},\ddot{u}) \right)$$

$$-D\left( \lambda^T f_{\dot{u}}(t,u,\dot{u},\ddot{u}) + y(t)^T g_{\dot{u}}(t,u,\dot{u},\ddot{u}) \right)$$

$$+D^2\left( \lambda^T f_{\ddot{u}}(t,u,\dot{u},\ddot{u}) + y(t)^T g_{\ddot{u}}(t,u,\dot{u},\ddot{u}) \right) = 0 \ ,\ t \in I,$$

$$\overline{z}^i(t)^T B^i(t) \overline{z}^i(t) \leqq 1 \ ,\ t \in I \ ,\ i \in P \ ,$$

$$\int_I y(t)^T g(t,u,\dot{u},\ddot{u}) dt \geqq 0 \ ,$$

$$y(t) \geqq 0 \ ,\ t \in I \ ,$$

$$\lambda > 0 \ .$$

Mond-Weir established usual duality theorems under the hypotheses that

$$\sum_{i=1}^{p} \lambda^i \int_I \left( f^i(t,.,.,.) + (\cdot)^T B^i(t) z^i(t) \right) dt$$

is pseudoinvex and $\int_I y(t)^T g(t,.,.,.) dt$ is quasi-invex with respect to the same $\eta$ .

In this paper, we propose, in the spirit of Husain and Jabeen [7] and Xu [18], a mixed type dual to (VP) and established various duality results under generalized invexity requirements. Special cases are generalized and it is also pointed out that our duality results can be viewed as dynamic generalizations of those of static case, already existing in the literature.

## 2. PRE-REQUISITES

Consider the real interval $I=[a,b]$, and the continuously differentiable function $\phi: I\times R^n\times R^n\times R^n\to R$. In order to consider $\phi(t,x,\dot{x},\ddot{x})$, where $x: I\to R^n$ is twice differentiable with its first and second order derivatives $\dot{x}$ and $\ddot{x}$ respectively, we denote the partial derivative $\phi$ by $\phi_t$ .

$$\phi_x=\left[\frac{\partial\phi}{\partial x^1},\ldots,\frac{\partial\phi}{\partial x^n}\right]^T \ , \ \phi_{\dot{x}}=\left[\frac{\partial\phi}{\partial \dot{x}^1},\ldots,\frac{\partial\phi}{\partial \dot{x}^n}\right]^T \ , \phi_{\ddot{x}}=\left[\frac{\partial\phi}{\partial \ddot{x}^1},\ldots,\frac{\partial\phi}{\partial \ddot{x}^n}\right]^T .$$

The partial derivative of other function will be written similarly. Let K designates the space of piecewise functions $x: I\to R^n$ possessing derivatives $\dot{x}$ and $\ddot{x}$ with the norm $\|x\|=\|x\|_\infty+\|Dx\|_\infty+\|D^2x\|_\infty$, where the differentiation operator $D$ is given by

$$u=Dx\Leftrightarrow x(t)=\alpha+\int_a^t u(s)\,ds,$$

Where $\alpha$ is given boundary value; thus $D\equiv\frac{d}{dt}$ except at discontinuities.

In the results to follow, we use $C(I,R^k)$ to denote the space of continuous functions $\phi: I\to R^m$ with the uniform norm; $T$ denotes matrix transpose.

Before stating our mixed type multiobjective variational problem, we mention the following conventions for vectors $x$ and $y$ in *n*-dimensional Euclidian space $R^n$ to be used throughout the analysis of this research together with some definitions of invexity and generalized invexity for easy reference.

$$x<y, \quad \Leftrightarrow \quad x_i<y_i, \quad i=1,2,\ldots,n.$$

$$x\leqq y, \quad \Leftrightarrow \quad x_i\leq y_i, \quad i=1,2,\ldots,n.$$

$$x\leq y, \quad \Leftrightarrow \quad x_i\leqq y_i, \quad i=1,2,\ldots,n, \text{ but } x\neq y$$

$$x\nleq y, \text{ is the negation of } x\leq y$$

For $x,y\in R$ , $x\leq y$ and $x<y$ have the usual meaning.

**DEFINITION 1.(*Invexity*)**: If there exists vector function $\eta(t,x,u)\in R^n$ with $\eta=0$ and $x(t)=u(t), t\in I$ and $D\eta=0$ for $\dot{x}(t)=\dot{u}(t), t\in I$ such that for a scalar function $\phi(t,x,\dot{x},\ddot{x})$ , the functional $\Phi(x,\dot{x},\ddot{x})=\int_I\phi(t,x,\dot{x},\ddot{x})\,dt$ satisfies

$$\Phi(x,\dot{u},\ddot{u}) - \Phi(x,\dot{x},\ddot{x}) \geqq \int_I \left\{ \eta^T \phi_x(t,x,\dot{x},\ddot{x}) + (D\eta)^T \phi_{\dot{x}}(t,x,\dot{x},\ddot{x}) + \left(D^2\eta\right)^T \phi_{\ddot{x}}(t,x,\dot{x},\ddot{x}) \right\} dt,$$

$\Phi$ is said to be invex in $x, \dot{x}$ *and* $\ddot{x}$ on $I$ with respect to $\eta$.

**DEFINITION 2. (***Pseudoinvexity***):** $\Phi$ is said to be pseudoinvex in $x, \dot{x}$ *and* $\ddot{x}$ with respect to $\eta$ if

$$\int_I \left\{ \eta^T \phi_x(t,x,\dot{x},\ddot{x}) + (D\eta)^T \phi_{\dot{x}}(t,x,\dot{x},\ddot{x}) + \left(D^2\eta\right)^T \phi_{\ddot{x}}(t,x,\dot{x},\ddot{x}) \right\} dt \geqq 0$$

implies $\quad \Phi(x,\dot{u},\ddot{u}) \geqq \Phi(x,\dot{x},\ddot{x})$.

**DEFINITION 3. (***Quasi-invex***):** The functional $\Phi$ is said to quasi-invex in $x, \dot{x}$ *and* $\ddot{x}$ with respect to $\eta$ if

$\Phi(x,\dot{u},\ddot{u}) \leqq \Phi(x,\dot{x},\ddot{x})$ implies

$$\int_I \left\{ \eta^T \phi_x(t,x,\dot{x},\ddot{x}) + (D\eta)^T \phi_{\dot{x}}(t,x,\dot{x},\ddot{x}) + \left(D^2\eta\right)^T \phi_{\ddot{x}}(t,x,\dot{x},\ddot{x}) \right\} dt \leqq 0.$$

We require the following definition of efficient solution for our further analysis. The set of feasible solutions for (VP) is designated by $X$.

**DEFINITION 4: (***Efficiency***)** A point $\bar{x} \in X$ is said to be efficient solution of (VP) if for all $\bar{x} \in X$

$$\int_I \left( f^i\left(t, x(t), \dot{x}(t), \ddot{x}(t)\right) dt + \left(x(t)^T B^i(t) x(t)\right)^{\frac{1}{2}} \right) dt \not\leq \int_I \left( f^i\left(t, \bar{x}(t), \dot{\bar{x}}(t), \ddot{\bar{x}}(t)\right) dt + \left(\bar{x}(t)^T B^i(t) \bar{x}(t)\right)^{\frac{1}{2}} \right) dt$$

for all $i \in P$.

In order to prove the strong duality theorem we will invoke the following lemma due to Changkong and Haimes [5].

**Lemma 1.** A point $\bar{x} \in X$ be an efficient solution of (VP) if and only if $\bar{x} \in X$ is an optimal solution of the following problem $\left(\mathrm{P}_k(\bar{x})\right)$ for all $k$.

$$\left(\mathrm{P}_k(\bar{x})\right) \quad \text{Minimize} \left( \int_I \left( f^k(t,x,\dot{x},\ddot{x}) dt + \left(x(t)^T B^k(t) x(t)\right)^{\frac{1}{2}} \right) dt \right.$$

Subject to

$$x(a) = 0 = x(b),$$

$$\dot{x}(a)=0=\dot{x}(b),$$

$$g(t,x,\dot{x},\ddot{x}) \leqq 0 \ , \ t\in I,$$

$$\int_I \left( f^i\left(t,x(t),\dot{x}(t),\ddot{x}(t)\right)dt + \left(x(t)^T B^i(t)\, x(t)\right)^{\frac{1}{2}} \right)dt$$

$$\leq \int_I \left( f^i\left(t,\bar{x}(t),\dot{\bar{x}}(t),\ddot{\bar{x}}(t)\right)dt + \left(\bar{x}(t)^T B^i(t)\,\bar{x}(t)\right)^{\frac{1}{2}} \right)dt \ , \ i\neq k$$

## 2. MIXED TYPE DUALITY

We formulate the following type dual (Mix D) to (VP):

**(Mix D)** Maximize $\left( \int_I \left( f^1(t,u,\dot{u},\ddot{u}) + u(t)^T B^1(t)\, z^1(t) + \sum_{j\in J_\circ} y^j(t)\, g^j(t,u,\dot{u},\ddot{u}) \right)dt \right.$

$$\left. ,\ldots, \int_I \left( f^p(t,u,\dot{u},\ddot{u}) + \left(u(t)^T B^p(t)\, z^p(t)\right) + \sum_{j\in J_\circ} y^j(t)\, g^j(t,u,\dot{u},\ddot{u}) \right)dt \right)$$

Subject to

$$u(a)=0=u(b), \qquad (1)$$

$$\dot{u}(a)=0=\dot{u}(b), \qquad (2)$$

$$\sum_{i=1}^{p} \lambda^i \left( f_u^i(t,u,\dot{u},\ddot{u})dt + B^i(t)\, z^i(t) + y(t)^T g_u(t,u,\dot{u},\ddot{u}) \right)$$

$$-D\left(\lambda^T f_{\dot{u}} + y(t)^T g_{\dot{u}}\right) + D^2\left(\lambda^T f_{\ddot{u}} + y(t)^T g_{\ddot{u}}\right) = 0 \ , \ t\in I \ , \qquad (3)$$

$$\sum_{j\in J_\alpha} \int_I y^j(t)\, g^j(t,u,\dot{u},\ddot{u})dt \geqq 0 \ , \ \alpha = 1,2,\ldots,r \, , \qquad (4)$$

$$\bar{z}^i(t)^T B^i(t)\, \bar{z}^i(t) \leqq 1 \ , \ t\in I \ , \ i\in P \ . \qquad (5)$$

$$y(t) \geqq 0 \ , \ t\in I \ , \qquad (6)$$

$$\lambda \in \Lambda^+ . \qquad (7)$$

where

(i) $\Lambda^+ = \left\{ \lambda\in R^p \,\middle|\, \lambda>0,\ \lambda^T e = 1, e = (1,1,\ldots,1)^T \in R^p \right\}$

(ii) $J_\alpha \subseteq M = \{1,2,\ldots,m\}$, $\alpha = 0,1,2,\ldots,r$ with $\bigcup_{\alpha=0}^{r} J_\alpha = M$ and $J_\alpha \bigcap J_\beta = \phi$ ,if $\alpha \neq \beta$.

If $J_\circ = M$ , then (Mix D) becomes Wolfe type dual considered in [9], if $J_\circ = \phi$ and $J_\alpha = M$ for some $\alpha \in \{1,2,\ldots,r\}$, then (Mix D) becomes Mond-Weir type dual considered in [9].

**THEOREM 1. (*Weak Duality*):** Let $\bar{x}$ be feasible for (VP) and $\left(u, y, z^1,\ldots,z^p,\lambda\right)$ be feasible for (Mix D). If for feasible $\left(x,u,z^1,\ldots,z^p,\lambda\right)$ ,

$\sum_{i=1}^{p} \lambda^i \int_I \left( f^i(t,.,.,.) + \sum_{j\in J_\circ} y^j(t) g^j(t,.,.,.) + (\cdot)^T B^i(t) z^i(t) \right) dt$ is pseudoinvex and

$\sum_{j\in J_\alpha} \int_I y^j(t) g^j(t,.,.,.)\, dt$ , $\alpha = 1,2,\ldots,r$ is quasi-invex with respect to same $\eta$, then the following cannot hold:

$$\int_I \left( f^i(t,x,\dot{x},\ddot{x})dt + \left( x(t)^T B^i(t) x(t) \right)^{\frac{1}{2}} \right) dt$$

$$\leqq \int_I \left( f^i(t,u,\dot{u},\ddot{u})dt + \left( u(t)^T B^i(t) z^i(t) \right) + \sum_{j\in J_\circ} y^j(t) g^j(t,u,\dot{u},\ddot{u}) \right) dt \tag{8}$$

For all $i \in P$ , and

$$\int_I \left( f^k(t,x,\dot{x},\ddot{x})dt + \left( x(t)^T B^k(t) x(t) \right)^{\frac{1}{2}} \right) dt$$

$$\leqq \int_I \left( f^k(t,u,\dot{u},\ddot{u})dt + \left( u(t)^T B^k(t) z^k(t) \right) + \sum_{j\in J_\circ} y^j(t) g^j(t,u,\dot{u},\ddot{u}) \right) dt \tag{9}$$

for some $k$ .

**PROOF:** Suppose contrary to the result, that (8) and (9) hold. In view of $y(t) \geqq 0$ , $t \in I$ , $\bar{z}^i(t)^T B^i(t) \bar{z}^i(t) \leqq 1$ , $t \in I$ , $i \in P$ and $g(t,x,\dot{x},\dot{x}) \leqq 0$ , $t \in I$,, these inequalities yield,

$$\int_I \left( f^i(t,x,\dot{x},\ddot{x})dt + \left( x(t)^T B^i(t) z^i(t) \right) + \sum_{j\in J_\circ} y^j(t) g^j(t,u,\dot{u},\ddot{u}) \right) dt$$

$$\leqq \int_I \left( f^i(t,u,\dot{u},\ddot{u})dt + \left( u(t)^T B^i(t) z^i(t) \right) + \sum_{j\in J_\circ} y^j(t) g^j(t,u,\dot{u},\ddot{u}) \right) dt$$

for all $i \in P$ , and

$$\int_I \left( f^k(t,x,\dot{x},\ddot{x})\,dt + \left(x(t)^T B^k(t) z^k(t)\right) + \sum_{j\in J_\circ} y^j(t) g^j(t,u,\dot{u},\ddot{u}) \right) dt$$

$$\leqq \int_I \left( f^k(t,u,\dot{u},\ddot{u})\,dt + \left(u(t)^T B^k(t) z^k(t)\right) + \sum_{j\in J_\circ} y^j(t) g^j(t,u,\dot{u},\ddot{u}) \right) dt$$

for some $k$ .

Now using $\lambda > 0$ and $\lambda^T e = 1$, these inequalities yield

$$\sum_{i=1}^{p} \lambda^i \int_I \left( f^i(t,x,\dot{x},\ddot{x})\,dt + \left(x(t)^T B^i(t) z^i(t)\right) + \sum_{j\in J_\circ} y^j(t) g^j(t,x,\dot{x},\ddot{x}) \right) dt$$

$$< \sum_{i=1}^{p} \lambda^i \int_I \left( f^i(t,u,\dot{u},\ddot{u})\,dt + \left(u(t)^T B^i(t) z^i(t)\right) + \sum_{j\in J_\circ} y^j(t) g^j(t,x,\dot{x},\ddot{x}) \right) dt$$

By pseudo invexity of $\sum_{i=1}^{p} \lambda^i \int_I \left( f^i(t,.,.,.) + \sum_{j\in J_\circ} y^j(t) g^j(t,.,.,.) + (\cdot)^T B^i(t) z^i(t) \right) dt$ with respect to $\eta$ , we have

$$0 > \sum_{i=1}^{p} \lambda^i \int_I \eta^T \left[ \left( f_x^i + B^i(t) z^i(t) + \sum_{j\in J_\circ} y^j(t) g_x^j \right) + (D\eta)^T \left( f_{\dot{x}}^i + \sum_{j\in J_\circ} y^j(t) g_{\dot{x}}^j \right) + (D^2\eta)^T \left( f_{\ddot{x}}^i + \sum_{j\in J_\circ} y^j(t) g_{\ddot{x}}^j \right) \right] dt$$

This, by integration by parts, gives,

$$= \sum_{i=1}^{p} \lambda^i \int_I \eta^T \left[ \left\{ \left( f_x^i + B^i(t) z^i(t) + \sum_{j\in J_\circ} y^j(t) g_x^j \right) - D\left( f_{\dot{x}}^i + \sum_{j\in J_\circ} y^j(t) g_{\dot{x}}^j \right) \right\} dt + \eta^T \left( f_{\dot{x}}^i + \sum_{j\in J_\circ} y^j(t) g_{\dot{x}}^j \right) \Bigg|_{t=a}^{t=b} + (D\eta)^T \left( f_{\ddot{x}}^i + \sum_{j\in J_\circ} y^j(t) g_{\ddot{x}}^j \right) \Bigg|_{t=a}^{t=b} + \int_I (D\eta)^T D\left( f_{\ddot{x}}^i + \sum_{j\in J_\circ} y^j(t) g_{\ddot{x}}^j \right) dt \right]$$

Using the boundary conditions which at $t = a$ , $t = b$ gives $D\eta = 0 = \eta$

$$= \sum_{i=1}^{p} \lambda^i \int_I \left[ \eta^T \left\{ \left( f_x^i + B^i(t) z^i(t) + \sum_{j\in J_\circ} y^j(t) g_x^j \right) - D\left( f_{\dot{x}}^i + \sum_{j\in J_\circ} y^j(t) g_{\dot{x}}^j \right) \right\} dt + \int_I (D\eta)^T D\left( f_{\ddot{x}}^i + \sum_{j\in J_\circ} y^j(t) g_{\ddot{x}}^j \right) \right] dt$$

$$= \sum_{i=1}^{p} \lambda^i \int_I \left[ \eta^T \left\{ \left( f_x^i + B^i(t) z^i(t) + \sum_{j \in J_\circ} y^j(t) g_x^j \right) - D\left( f_{\dot{x}}^i + \sum_{j \in J_\circ} y^j(t) g_{\dot{x}}^j \right) \right.\right.$$

$$\left.\left. + D^2 \left( f_{\ddot{x}}^i + \sum_{j \in J_\circ} y^j(t) g_{\ddot{x}}^j \right) \right\} \right] dt - \eta^T \left( f_{\ddot{x}}^i + \sum_{j \in J_\circ} y^j(t) g_{\ddot{x}}^j \right) \Bigg|_{t=a}^{t=b}$$

$$0 > \sum_{i=1}^{p} \lambda^i \int_I \left[ \eta^T \left\{ \left( f_x^i + B^i(t) z^i(t) + \sum_{j \in J_\circ} y^j(t) g_x^j \right) - D\left( f_{\dot{x}}^i + \sum_{j \in J_\circ} y^j(t) g_{\dot{x}}^j \right) \right.\right.$$

$$\left.\left. + D^2 \left( f_{\ddot{x}}^i + \sum_{j \in J_\circ} y^j(t) g_{\ddot{x}}^j \right) \right\} \right] dt \tag{10}$$

Now, from the feasibility of (VP) and (Mix-D), we have

$$\sum_{j \in J_\alpha} \int_I y^j(t) g^j(t, x, \dot{x}, \ddot{x})\, dt \leqq \sum_{j \in J_\alpha} \int_I y^j(t) g^j(t, u, \dot{u}, \ddot{u})\, dt$$

This, because of quasi-invexity of $\sum_{j \in J_\alpha} \int_I y^j(t) g^j(t,.,.,.)\, dt$ with $\eta$ implies

$$\sum_{j \in J_\alpha} \int_I \left\{ \eta^T \left( y^j(t) g_x^j \right) + (D\eta)^T y^j(t) g_{\dot{x}}^j + \left( D^2 \eta \right)^T y^j(t) g_{\ddot{x}}^j \right\} dt \leqq 0$$

As earlier, integrating by parts and using the boundary conditions, we have

$$\sum_{j \in J_\alpha} \int_I \eta^T \left\{ \left( y^j(t) g_x^j \right) + D^T y^j(t) g_{\dot{x}}^j + D^2 y^j(t) g_{\ddot{x}}^j \right\} dt \leqq 0 \tag{11}$$

Combining (10) and (11), we have

$$\int_I \eta^T \left( \sum_{i=1}^{p} \lambda^i \left( f_x^i + B^i(t) z^i(t) + y(t)^T g_x \right) \right.$$

$$\left. - D\left( \lambda^T f_{\dot{x}} + y(t)^T g_{\dot{x}} \right) + D^2 \left( \lambda^T f_{\ddot{x}} + y(t)^T g_{\ddot{x}} \right) \right) dt < 0$$

From the equality constraint of the dual, we have

$$\int_I \eta^T \left( \sum_{i=1}^{p} \lambda^i \left( f_x^i + B^i(t) z^i(t) + y(t)^T g_x \right) \right.$$

$$\left. - D\left( \lambda^T f_{\dot{x}} + y(t)^T g_{\dot{x}} \right) + D^2 \left( \lambda^T f_{\ddot{x}} + y(t)^T g_{\ddot{x}} \right) \right) dt = 0$$

which is a contradiction. Hence the conclusion of the theorem is true.

**THEOREM 2.** (*Strong Duality*)**:** Let $\bar{x} \in X$ be an efficient solution of (VP) and for at least one $i \in P$, $\bar{x}$ satisfies the regularity condition [3] for the problem $\left( \mathrm{P}_k(\bar{x}) \right)$. Then there exist multipliers $\bar{\lambda} \in R^p$, piecewise smooth $\bar{y} \in R^m$ and $z^i(t) \in R^n$, $i = \{1, 2, \ldots, p\}$ such that $(\bar{x}, \bar{y}, \bar{z}^1, \ldots, \bar{z}^p, \lambda)$ is feasible for (Mix D) and the objectives of (VP) and (Mix D) are equal.

Further, if the hypotheses of Theorem 1 are met, then $(\bar{x}, \bar{y}, \bar{z}^1, \ldots, \bar{z}^p, \lambda)$ is an efficient solution of (Mix D).

**PROOF:** Since $\bar{x} \in X$ is an optimal solution of $(\mathrm{P}_k(\bar{x}))$.This implies that there exist $\bar{\xi} \in R^p$ with $\bar{\xi}^1, \ldots, \bar{\xi}^p, z^i(t) \in R^n$, $i = \{1, 2, \ldots, p\}$ and piecewise smooth $\bar{\upsilon} \in R^m$ such that, the following optimality conditions [3] hold:

$$\bar{\xi}^k \left( f_x^k + B^k(t) z^k(t) - D f_{\dot{x}}^k + D^2 f_{\ddot{x}}^k \right)$$

$$+ \sum_{\substack{i=1 \\ i \neq k}}^{p} \bar{\xi}^i \left( f_x^i + B^i(t) z^i(t) - D f_{\dot{x}}^i + D^2 f_{\ddot{x}}^i \right)$$

$$+ \bar{\upsilon}(t)^T g_x - D\left(\bar{\upsilon}(t)^T g_{\dot{x}}\right) + D^2\left(\bar{\upsilon}(t)^T g_{\ddot{x}}\right) = 0 \tag{12}$$

$$\bar{\upsilon}(t)^T g\left(t, \bar{x}, \dot{\bar{x}}, \ddot{\bar{x}}\right) dt = 0 \tag{13}$$

$$\left(\bar{x}(t)^T B^i(t) \bar{x}(t)\right)^{\frac{1}{2}} = \left(\bar{x}(t)^T B^i(t) \bar{z}^i(t)\right), \; i = 1, \ldots, p \tag{14}$$

$$\left(\bar{z}(t)^T B^i(t) \bar{z}^i(t)\right) \leqq 1, \; t \in I, \; i = 1, 2, \ldots, p$$

(15)

$$\bar{\xi} > 0 \;, \; \bar{\upsilon}(t) \geqq 0, \; t \in I \tag{16}$$

Dividing (12), (13) and (16) by $\sum_{i=1}^{p} \bar{\xi}^i$, and setting $\bar{\lambda}^i = \dfrac{\bar{\xi}}{\sum_{i=1}^{p} \bar{\xi}^i}$, $i = 1, \ldots, p$ and $\bar{y}(t) = \dfrac{\bar{\upsilon}(t)}{\sum_{i=1}^{p} \bar{\xi}^i}$, we have

$$\sum_{i=1}^{p} \bar{\lambda}^i \left( f_x^i + B^i(t) \bar{z}^i(t) - D f_{\dot{x}}^i + D^2 f_{\ddot{x}}^i \right) + \bar{y}(t)^T g_x - D\left(\bar{y}(t)^T g_{\dot{x}}\right) + D^2\left(\bar{y}(t)^T g_{\ddot{x}}\right) = 0 \tag{17}$$

$$\bar{y}(t)^T g\left(t, \bar{x}, \dot{\bar{x}}, \ddot{\bar{x}}\right) dt = 0 \tag{18}$$

$$\bar{\lambda} > 0 \;, \; \bar{y}(t) \geqq 0, \; t \in I \tag{19}$$

The equation (18) implies

$$\sum_{j\in J_\circ} y^j(t) g^j(t,x,\dot{x},\ddot{x}) = 0 \text{ and } \sum_{j\in J_\alpha} y^j(t) g^j(t,x,\dot{x},\ddot{x}) = 0 \;,\; \alpha = 1,\dots,r \tag{20}$$

This implies,

$$\sum_{j\in J_\alpha} \int_I y^j(t) g^j(t,x,\dot{x},\ddot{x}) = 0 \;,\; \alpha = 1,\dots,r \tag{21}$$

Consequently (15), (17) ,(19) and (21) implies that $\left(\bar{x},\bar{y},\bar{z}^1,\dots,\bar{z}^p,\bar{\lambda}\right)$ is feasible for (Mix D).

$$\int_I \left( f^i\left(t,\bar{x},\dot{\bar{x}},\ddot{\bar{x}}\right) + \left(\bar{x}(t)^T B^i(t)\bar{x}(t)\right)^{\frac{1}{2}} + \sum_{j\in J_\circ} y^j(t) g^j\left(t,\bar{x},\dot{\bar{x}},\ddot{\bar{x}}\right) \right) dt$$

$$= \int_I \left( f^i\left(t,\bar{x},\dot{\bar{x}},\ddot{\bar{x}}\right) + \left(\bar{x}(t)^T B^i(t)\bar{z}^i(t)\right)\right) dt \;,\; i = 1,2,\dots,p$$

This in view of Theorem 1, the efficiency of $\left(\bar{x},\bar{y},\bar{z}^1,\dots,\bar{z}^p,\bar{\lambda}\right)$ follows.

As in [15], by employing chain rule in calculus, it can be easily seen that the expression

$$\sum_{i=1}^{p} \lambda^i \left( f_x^i(t,u,\dot{u},\ddot{u}) dt + B^i(t) z^i(t) + y(t)^T g_x(t,u,\dot{u},\ddot{u}) \right)$$

$$-D\left(\lambda^T f_{\dot{x}} + y(t)^T g_{\dot{x}}\right) + D^2\left(\lambda^T f_{\ddot{x}} + y(t)^T g_{\ddot{x}}\right) = 0, \quad , t \in I$$ may be regarded as a function $\theta$ of variables $t, x, \dot{x}, \ddot{x}, \dddot{x}, y, \dot{y}, \ddot{y}$ and $\lambda$ , where $\dddot{x} = D^3 x$ and $\ddot{y} = D^2 y$ . That is, we can write

$$\theta(t,x,\dot{x},\ddot{x},\dddot{x},y,\dot{y},\ddot{y},\lambda) = \sum_{i=1}^{p} \lambda^i \left( f_x^i(t,u,\dot{u},\ddot{u}) dt + B^i(t) z^i(t) + y(t)^T g_x(t,u,\dot{u},\ddot{u}) \right)$$

$$- D\left(\lambda^T f_{\dot{x}} + y(t)^T g_{\dot{x}}\right) + D^2\left(\lambda^T f_{\ddot{x}} + y(t)^T g_{\ddot{x}}\right) = 0$$

In order to prove converse duality between (VP) and (Mix D), the space $X$ is now replaced by a smaller space $X_2$ of piecewise smooth thrice differentiable function $x : I \to R^n$ with the norm $\|x\|_\infty + \|Dx\|_\infty + \|D^2 x\|_\infty + \|D^3 x\|_\infty$ .The problem (Mix D) may now be briefly written as,

$$\text{Minimize}\left( -\int_I \left( f^1(t,u,\dot{u},\ddot{u}) + u(t)^T B^1(t) z^1(t) + \sum_{j\in J_\circ} y^j(t) g^j(t,u,\dot{u},\ddot{u}) \right) dt \right.$$

$$\left. ,\dots, \int_I -\left( f^p(t,u,\dot{u},\ddot{u}) + \left(u(t)^T B^p(t) z^p(t)\right) + \sum_{j\in J_\circ} y^j(t) g^j(t,u,\dot{u},\ddot{u}) \right) dt \right)$$

Subject to

$$u(a)=0=u(b)$$

$$\dot{u}(a)=0=\dot{u}(b)$$

$$\theta(t,x,\dot{x},\ddot{x},\dddot{x},y,\dot{y},\ddot{y},\lambda)=0\ ,\ t\in I$$

$$\sum_{j\in J_\alpha}\int_I y^j(t)g^j(t,u,\dot{u},\ddot{u})dt \geqq 0\ ,\ \alpha=1,2,\ldots,r$$

$$\bar{z}^i(t)^T B^i(t)\bar{z}^i(t) \leqq 1\ ,\ t\in I\ ,\ i\in P$$

$$y(t)\geqq 0\ ,\ t\in I$$

$$\lambda>0\ ,\ \lambda^T e=1$$

where

$$\theta(t,x,\dot{x},\ddot{x},\dddot{x},y,\dot{y},\ddot{y},\lambda)=\sum_{i=1}^{p}\lambda^i\left(f_x^i(t,u,\dot{u},\ddot{u})dt+B^i(t)z^i(t)+y(t)^T g_x(t,u,\dot{u},\ddot{u})\right)$$
$$-D\left(\lambda^T f_{\dot{x}}+y(t)^T g_{\dot{x}}\right)+D^2\left(\lambda^T f_{\ddot{x}}+y(t)^T g_{\ddot{x}}\right)=0$$

Consider $\theta(t,x(\cdot),\dot{x}(\cdot),\ddot{x}(\cdot),\dddot{x}(\cdot),y(\cdot),\dot{y}(\cdot),\ddot{y}(\cdot),\lambda)=0$ as defining a mapping $\psi: X\times Y\times R^p\to B$ where $Y$ is a space of piecewise twice differentiable function and $B$ is the Banach Space. In order to apply Theorem 1[10] to the problem (Mix D), the infinite dimensional inequality must be restricted. In the following theorem, we use $\psi'$ to represent the Frèchèt derivative $\left[\psi_x(x,y,\lambda),\psi_y(x,y,\lambda),\psi_\lambda(x,y,\lambda)\right]$.

**THEOREM 3. (*Converse Duality*):** Let $\left(\bar{x},\bar{y},\bar{z}^1,\ldots,\bar{z}^p,\bar{\lambda}\right)$ be an efficient solution for (Mix D). Assume that

**(A$_1$)** The Frèchèt derivative $\psi'$ has a (weak*) closed range,

**(A$_2$)** $f$ and $g$ are twice continuously differentiable,

**(A$_3$)** $$\left(f_x^i+B^i(t)z^i(t)+\sum_{j\in J_\circ}y^j(t)g_x^i\right)-D\left(f_{\dot{x}}^i+\sum_{j\in J_\circ}y^j(t)g_{\dot{x}}^i\right)$$
$$+D^2\left(f_{\ddot{x}}^i+\sum_{j\in J_\circ}y^j(t)g_{\ddot{x}}^i\right),\ i\in\{1,2,\ldots,p\}$$ are linearly independent, and

**(A$_4$)** $$\left(\beta(t)^T\theta_x-D\beta(t)^T\theta_{\dot{x}}+D^2\beta(t)^T\theta_{\ddot{x}}-D^3\beta(t)^T\theta_{\dddot{x}}\right)\beta(t)=0,\Rightarrow\beta(t)=0\ ,\ t\in I$$

Further, if the hypotheses of Theorem 1 are met then $\bar{x}$ is an efficient solution of (VP).

**PROOF:** Since $\left(\bar{x},\bar{y},\bar{z}^1,\ldots,\bar{z}^p,\bar{\lambda}\right)$ with $\psi'$ having a (weak*) closed range, is an efficient solution of (Mix D), then by Theorem 1[9], there exist multipliers $\tau \in R^p$ , piecewise smooth $\beta(t): I \to R^n$ , $z^i(t) \in R^n$ , $i=1,\ldots,p$, $\gamma \in R$ for each of $r$ constraints, $\eta \in R^p$, $\delta \in R$ and piecewise smooth $\mu(t): I \to R^m$ satisfying the following Fritz-John optimality conditions [9].

$$-\sum_{i=1}^{p} \tau^i \left[ \left( f_x^i + B^i(t) z^i(t) + \sum_{j\in J_\circ} y^j(t) g_x^j \right) - D\left( f_{\dot{x}}^i + \sum_{j\in J_\circ} y^j(t) g_{\dot{x}}^j \right) + D^2\left( f_{\ddot{x}}^i + \sum_{j\in J_\circ} y^j(t) g_{\ddot{x}}^j \right) \right] - \gamma \sum_{\alpha=1}^{r} \sum_{j\in J_\alpha} \left( y^j(t) g_x^j - D y^j(t) g_{\dot{x}}^j + D^2 y^j(t) g_{\ddot{x}}^j \right) + \beta(t)^T \theta_x - D\beta(t)^T \theta_{\dot{x}} + D^2 \beta(t)^T \theta_{\ddot{x}} - D^3 \beta(t)^T \theta_{\dddot{x}} = 0 ,\ t \in I \tag{22}$$

$$-\left(\tau^T e\right) g^j + \beta(t)^T \theta_{y^j} - D\beta(t)^T \theta_{\dot{y}^j} + D^2 \beta(t)^T \theta_{\ddot{y}^j} - \mu^j(t) = 0 \ ,\ j \in J_\circ$$

(23)

$$-\gamma g^j + \beta(t)^T \theta_{y^j} - D\beta(t)^T \theta_{\dot{y}^j} + D^2 \beta(t)^T \theta_{\ddot{y}^j} - \mu^j(t) = 0 \ ,\ j \in J_\alpha ,\alpha = 1,\ldots,r \tag{24}$$

$$-\tau^i x(t)^T B^i(t) + \lambda^i \beta(t)^T B^i(t) + 2\delta^i B^i(t) z^i(t) = 0 \ ,\ i = 1,\ldots,p \tag{25}$$

$$\left[ f_x^i + B^i(t) z^i(t) + \sum_{j\in J_\circ} y^j(t) g_x^j - D\left( f_{\dot{x}}^i + \sum_{j\in J_\circ} y^j(t) g_{\dot{x}}^j \right) + D^2\left( f_{\ddot{x}}^i + \sum_{j\in J_\circ} y^j(t) g_{\ddot{x}}^j \right) \right] - \eta^i = 0 \tag{26}$$

$$\mu^T(t)\bar{y}(t) = 0 \ ,\ t \in I \tag{27}$$

$$\eta^T \lambda = 0$$

(28)

$$\gamma \sum_{j\in J_\alpha} \int_I y(t)^T g\left(t,\bar{x},\dot{\bar{x}},\ddot{\bar{x}}\right) dt = 0 \ ,\ \alpha = 1,2,\ldots,p \tag{29}$$

$$\delta^i \left( z^i(t)^T B^i(t) z^i(t) - 1 \right) = 0 \ ,\ t \in I \tag{30}$$

$$\left(\tau,\gamma,\mu(t),\delta,\eta\right) \geqq 0 \tag{31}$$

$$\left(\tau,\beta(t),\gamma,\mu(t),\delta,\eta\right) \neq 0 \tag{32}$$

Since $\lambda > 0$, (28) implies $\eta = 0$. Consequently (26) implies

$$\left[ f_x^i + B^i(t)z^i(t) + \sum_{j \in J_\circ} y^j(t) g_x^j - D\left( f_{\dot{x}}^i + \sum_{j \in J_\circ} y^j(t) g_{\dot{x}}^j \right) + D^2\left( f_{\ddot{x}}^i + \sum_{j \in J_\circ} y^j(t) g_{\ddot{x}}^j \right) \right] = 0 \qquad (33)$$

Using the duality constraint of (Mix D) equation (22) reduces to

$$\begin{aligned} -\sum_{i=1}^{p}\left(\tau^i - \gamma\lambda^i\right)\Bigg\{ & \left( f_x^i + B^i(t)z^i(t) + \sum_{j \in J_\circ} y^j(t) g_x^j \right) - D\left( f_{\dot{x}}^i + \sum_{j \in J_\circ} y^j(t) g_{\dot{x}}^j \right) \\ & + D^2\left( f_{\ddot{x}}^i + \sum_{j \in J_\circ} y^j(t) g_{\ddot{x}}^j \right)\Bigg\} + \beta(t)^T \theta_x - D\beta(t)^T \theta_{\dot{x}} + D^2\beta(t)^T \theta_{\ddot{x}} - D^3\beta(t)^T \theta_{\dddot{x}} = 0 \text{ , } t \in I \end{aligned} \qquad (34)$$

Post multiplying (34) by $\beta(t)$ and then using (33), we obtain,

$$\left(\beta(t)^T \theta_x - D\beta(t)^T \theta_{\dot{x}} + D^2\beta(t)^T \theta_{\ddot{x}} - D^3\beta(t)^T \theta_{\dddot{x}}\right)\beta(t) = 0 \text{ , } t \in I$$

This because of the hypothesis ($A_4$) yields,

$$\beta(t) = 0 \text{ , } t \in I \qquad (35)$$

Using (35) in (34) in (22), we have

$$\begin{aligned} -\sum_{i=1}^{p}\left(\tau^i - \gamma\lambda^i\right)\Bigg\{ & \left( f_x^i + B^i(t)z^i(t) + \sum_{j \in J_\circ} y^j(t) g_x^j \right) - D\left( f_{\dot{x}}^i + \sum_{j \in J_\circ} y^j(t) g_{\dot{x}}^j \right) \\ & + D^2\left( f_{\ddot{x}}^i + \sum_{j \in J_\circ} y^j(t) g_{\ddot{x}}^j \right)\Bigg\} = 0 \end{aligned}$$

This, due to the hypothesis ($A_3$) gives,

$$\tau^i - \gamma\lambda^i = 0 \text{ , } i = 1, 2, \ldots, p \qquad (36)$$

Suppose $\gamma = 0$, then from (36) we have $\tau = 0$. Consequently from (23) and (24) implies $\mu(t) = 0$ , $t \in I$. Also from (25) we have $\delta^i B^i(t)z^i(t) = 0$, which together with (30) implies $\delta^i = 0$, i.e. $\delta = 0$. Thus, $\left(\tau, \beta(t), \mu(t), \eta, \gamma, \delta,\right) = 0$, which is a contradiction to (32). Hence $\gamma > 0$. From (34) it implies that $\tau > 0$.

By the generalized Schwartz inequality [17],

$$\left(\bar{x}(t)^T B^i(t)\bar{z}^i(t)\right) \leq \left(\bar{x}(t)^T B^i(t)\bar{x}(t)\right)^{\frac{1}{2}} \left(\bar{z}^i(t) B^i(t)\bar{z}^i(t)\right)^{\frac{1}{2}} \text{ , } i = 1, 2, \ldots, p \qquad (37)$$

Now let $\frac{\delta^i}{\tau^i} = \alpha^i$. Then $\alpha^i \geqq 0$ and from (25), we have

$$B^i(t)\bar{x}(t) = 2\alpha^i B^i(t)\bar{z}^i(t), \ i = 1, 2, \ldots, p$$

which, is the condition for the equality in (37). Therefore, we have

$$\left(\bar{x}(t)^T B^i(t)\bar{z}^i(t)\right) = \left(\bar{x}(t)^T B^i(t)\bar{x}(t)\right)^{\frac{1}{2}} \left(z^i(t) B^i(t) z^i(t)\right)^{\frac{1}{2}}$$

From (30), either $\delta^i = 0$ or $\bar{z}^i(t)^T B^i(t)\bar{z}^i(t) = 1$ or $\delta^i = 0$, $i = 1, \ldots, p$ and hence $B^i(t)\bar{x}(t) = 0.$ Therefore, in either case

$$\left(x(t)^T B^i(t) z^i(t)\right) = \left(x(t)^T B^i(t) x(t)\right)^{\frac{1}{2}}, \ i = 1, 2, \ldots, p \tag{38}$$

From (23) and (24) we readily obtain,

$$g\left(t, \bar{x}, \dot{\bar{x}}, \ddot{\bar{x}}\right) \leqq 0, \ t \in I$$

which gives the feasibility of $\bar{x}$ for (P). Using (27) in (23) and (24) we have

$$\bar{y}(t)^T g = 0, \ \text{ i.e., } \bar{y}^j(t) g^j = 0, \ j = 1, 2, \ldots, m, \ t \in I$$

This obviously gives

$$\sum_{j \in J_\circ} \bar{y}^j(t) g^j = 0, \ j = 1, 2, \ldots, m \tag{39}$$

Hence

$$\int_I \left( f^i\left(t, \bar{x}, \dot{\bar{x}}, \ddot{\bar{x}}\right) dt + \bar{x}(t)^T B^i(t)\bar{z}^i(t) + \sum_{j \in J_\circ} \bar{y}^j(t) g^j = 0 \right) dt$$

$$= \int_I \left( f^i\left(t, \bar{x}, \dot{\bar{x}}, \ddot{\bar{x}}\right) dt + \left(\bar{x}(t)^T B^i(t)\bar{x}(t)\right)^{\frac{1}{2}} \right) dt, \ i = 1, 2, \ldots, p$$

(by using (38) and (39))

The efficiency of $\bar{x}$ for (VP) follows by an application of Theorem 1.

## 4. VARIATIONAL PROBLEMS WITH NATURAL BOUNDARY CONDITIONS

It is possible to extend the duality theorems established in the previous sections to the corresponding variational problems with natural boundary values rather that fixed points.

**(VP)$_0$:** Minimize
$$\left( \int_I \left( f^1(t, x, \dot{x}, \ddot{x}) dt + \left(x(t)^T B^1(t) x(t)\right)^{\frac{1}{2}} \right) dt, \ldots, \int_I \left( f^p(t, x, \dot{x}, \ddot{x}) dt + \left(x(t)^T B^p(t) x(t)\right)^{\frac{1}{2}} \right) dt \right)$$

Subject to

$$g^{j}\left(t,x,\dot{x},\dot{x}\right)\leqq 0\ ,\ t\in I,\quad j=1,\ldots,m$$

**(MIX D)$_0$:** Maximize $\left(\int_I\left(f^{1}\left(t,u,\dot{u},\ddot{u}\right)+u\left(t\right)^{T}B^{1}\left(t\right)z^{1}\left(t\right)+\sum_{j\in J_\circ}y^{j}\left(t\right)g^{j}\left(t,u,\dot{u},\ddot{u}\right)\right)dt\right.$

$$\left.,\ldots,\int_I\left(f^{p}\left(t,u,\dot{u},\ddot{u}\right)+\left(u\left(t\right)^{T}B^{p}\left(t\right)z^{p}\left(t\right)\right)+\sum_{j\in J_\circ}y^{j}\left(t\right)g^{j}\left(t,u,\dot{u},\ddot{u}\right)\right)dt\right)$$

Subject to

$$\sum_{i=1}^{p}\lambda^{i}\left(f_u^{i}\left(t,u,\dot{u},\ddot{u}\right)+B^{i}\left(t\right)z^{i}\left(t\right)+y\left(t\right)^{T}g_u\left(t,u,\dot{u},\ddot{u}\right)\right)$$

$$-D\left(\lambda^{T}f_{\dot{u}}\left(t,u,\dot{u},\ddot{u}\right)+y\left(t\right)^{T}g_{\dot{u}}\left(t,u,\dot{u},\ddot{u}\right)\right)$$

$$+D^{2}\left(\lambda^{T}f_{\ddot{u}}\left(t,u,\dot{u},\ddot{u}\right)+y\left(t\right)^{T}g_{\ddot{u}}\left(t,u,\dot{u},\ddot{u}\right)\right)=0\ ,\ t\in I,$$

$$\left.\left(f_{\dot{x}}^{i}+\sum_{j\in J_\circ}y^{j}\left(t\right)g_{\dot{x}}^{j}\right)\right|_{t=a}=0\,,$$

$$\left.\left(f_{\dot{x}}^{i}+\sum_{j\in J_\circ}y^{j}\left(t\right)g_{\ddot{x}}^{j}\right)\right|_{t=b}=0\,,$$

$$\left.\left(f_{\ddot{x}}^{i}+\sum_{j\in J_\circ}y^{j}\left(t\right)g_{\dot{x}}^{j}\right)\right|_{t=a}=0\,,$$

$$\left.\left(f_{\ddot{x}}^{i}+\sum_{j\in J_\circ}y^{j}\left(t\right)g_{\ddot{x}}^{j}\right)\right|_{t=b}=0\,,$$

$$\forall\ i\in\left\{1,2,\ldots,p\right\}$$

$$\bar{z}^{i}\left(t\right)^{T}B^{i}\left(t\right)\bar{z}^{i}\left(t\right)\leqq 1\ ,\ t\in I\ ,\ i\in P\,,$$

$$\sum_{j\in J_\alpha}\int_I y^{j}\left(t\right)g^{j}\left(t,u,\dot{u},\ddot{u}\right)dt\geqq 0\ \ ,\ \ \alpha=1,2,\ldots,r\,,$$

$$\lambda>0\ \ ,\quad y\left(t\right)\geqq 0\ ,\ t\in I\ \ .$$

## 5. NONLINEAR PROGRAMMING.

If all the functions are independent of $t$ then the problems (VP)$_0$ and (MIX D)$_0$ become the following problems.

**(VP)$_1$:** Minimize$\left( f^1(x)+\left(x^T B^1 x\right)^{\frac{1}{2}},\ldots, f^p(x)+\left(x^T B^p x\right)^{\frac{1}{2}} \right)$

Subject to

$$g(x) \leqq 0$$

**(Mix D)$_1$:**

Maximize$\left( f^1(u)+u^T B^1 z^1 + \sum_{j\in J_\circ} y^j g^j(u),\ldots, f^P(u)+u^T B^P z^P + \sum_{j\in J_\circ} y^j g^j(u) \right)$

Subject to

$$\sum_{i=1}^{p} \lambda^i \left( f_u^i(u)dt + B^i z^i + y^T g_u(u) \right) = 0$$

$$\overline{z}^{iT} B^i \overline{z}^i \leqq 1 \quad , i \in P$$

$$\sum_{j\in J_\alpha} y^j g^j(u) \geqq 0 \quad , \ \alpha = 1,2,\ldots,r$$

$$y(t) \geqq 0 \ , t \in I$$

$$\lambda \in \Lambda^+$$

where $\quad \Lambda^+ = \left\{ \lambda \in R^p \,\middle|\, \lambda > 0,\ \lambda^T e = 1, e = (1,1,\ldots,1)^T \in R^p \right\}$

The duality results for this pair of problems are not explicitly reported in the literature but can be derived easily on the lines of the analysis of this research.